\documentclass[12pt,reqno]{amsart}
\usepackage{enumerate, latexsym, amsmath, amsfonts, amssymb, amsthm, color}
\def\pmod #1{\ ({\rm{mod}}\ #1)}
\def\Z{\Bbb Z}
\def\N{\Bbb N}

\def\l{\left}
\def\r{\right}
\def\bg{\bigg}
\def\({\bg(}
\def\){\bg)}
\def\t{\text}
\def\f{\frac}

\def\ls{\leqslant}
\def\gs{\geqslant}

\def\bi{\binom}

\def\eq{\equiv}

\def\Ack{\medskip\noindent {\bf Acknowledgment}}
\theoremstyle{plain}
\newtheorem{theorem}{Theorem}

\newtheorem{lemma}{Lemma}
\newtheorem{corollary}{Corollary}
\newtheorem{conjecture}{Conjecture}
\theoremstyle{definition}

\theoremstyle{remark}
\newtheorem{remark}{Remark}

 \vspace{4mm}

\begin{document}

\hbox{Int. J. Number Theory 12(2016), no.\,2, 527--539.}
\medskip

\title
[{Two congruences involving harmonic numbers}]
{Two congruences involving harmonic numbers with applications}

\author
[Guo-Shuai Mao and Zhi-Wei Sun] {Guo-Shuai Mao and Zhi-Wei Sun}

\address {(Guo-Shuai Mao) Department of Mathematics, Nanjing
University, Nanjing 210093, People's Republic of China}
\email{mg1421007@smail.nju.edu.cn}

\address{(Zhi-Wei Sun) Department of Mathematics, Nanjing
University, Nanjing 210093, People's Republic of China}
\email{zwsun@nju.edu.cn}

\keywords{Central binomial coefficients, congruences, harmonic numbers, Bernoulli polynomials.
\newline \indent 2010 {\it Mathematics Subject Classification}. Primary 11B65, 11B68; Secondary 05A10, 11A07.
\newline \indent The second author is the corresponding author. This research was supported by the Natural Science Foundation of China (grant 11571162).}

 \begin{abstract} The harmonic numbers $H_n=\sum_{0<k\ls n}1/k\ (n=0,1,2,\ldots)$ play important roles in mathematics. Let $p>3$ be a prime. With helps of some combinatorial identities, we establish the following two new congruences:
 $$\sum_{k=1}^{p-1}\f{\bi{2k}k}kH_k\eq\f13\l(\f p3\r)B_{p-2}\l(\f13\r)\pmod{p}$$
 and
 $$\sum_{k=1}^{p-1}\f{\bi{2k}k}kH_{2k}\eq\f7{12}\l(\f p3\r)B_{p-2}\l(\f13\r)\pmod{p},$$
 where $B_n(x)$ denotes the Bernoulli polynomial of degree $n$. As an application, we determine $\sum_{n=1}^{p-1}g_n$ and $\sum_{n=1}^{p-1}h_n$ modulo $p^3$, where
 $$g_n=\sum_{k=0}^n\bi nk^2\bi{2k}k\quad\t{and}\quad h_n=\sum_{k=0}^n\bi nk^2C_k$$
 with $C_k=\bi{2k}k/(k+1)$.
\end{abstract}

\maketitle

\section{Introduction}
\setcounter{lemma}{0}
\setcounter{theorem}{0}
\setcounter{corollary}{0}
\setcounter{remark}{0}
\setcounter{equation}{0}
\setcounter{conjecture}{0}

 For $n\in\N=\{0,1,2,\ldots\}$, define
 $$H_n:=\sum_{0<k\ls n}\f1k\quad\t{and}\quad H_n^{(2)}:=\sum_{0<k\ls n}\f1{k^2}.$$
 Those $H_n$ with $n\in\N$ are the classical harmonic numbers, and those $H_n^{(2)}$ with $n\in\N$
 are called the second-order harmonic numbers.

Let $p>3$ be a prime. By a classical result of J. Wolstenholme \cite{W}, we have
$$H_{p-1}\eq0\pmod{p^2}\ \t{and}\ H_{p-1}^{(2)}\eq0\pmod p,$$
which imply that
$$\f12\bi{2p}p=\bi{2p-1}{p-1}\eq1\pmod{p^3}.$$
Z.-W. Sun \cite{S12a} established some fundamental congruences involving harmonic numbers; for example, he showed that
$\sum_{k=1}^{p-1}H_k/(k2^k)\eq0\pmod p$ motivated by the known identity $\sum_{k=1}^\infty H_k/(k2^k)=\pi^2/12$.

Let $p>3$ be a prime. By Sun and R. Tauraso \cite[(1.9)]{ST11}, and Sun \cite[(2.9)]{S12c}, we have
$$\sum_{k=0}^{p-1}\bi{2k}k\eq\l(\f p3\r)\pmod{p^2}$$
and
$$\sum_{k=0}^{p-1}(-1)^k\bi{p-1}k\bi{2k}k\eq3^{p-1}\l(\f p3\r)\pmod{p^2}$$
respectively, where $(-)$ denotes the Legendre symbol. Hence
$$\sum_{k=0}^{p-1}\bi{2k}kH_k\eq\l(\f p3\r)\f{1-3^{p-1}}p\pmod p$$
since
$$(-1)^k\bi{p-1}k=\prod_{0<j\ls k}\l(1-\f pj\r)\eq1-pH_k\pmod{p^2}$$
for all $k=0,1,\ldots,p-1$.
In 2010 Sun and Tauraso \cite{ST10} proved that
$$\sum_{k=1}^{p-1}\f{\bi{2k}k}k\eq\f 89p^2B_{p-3}\pmod{p^3},$$
where $B_0,B_1,B_2,\ldots$ are the Bernoulli numbers given by
$$\f x{e^x-1}=\sum_{n=0}^\infty B_n\f{x^n}{n!}\ \ (0<|x|<2\pi).$$
In 2011 Sun \cite{S11} showed that
$$\sum_{k=1}^{(p-1)/2}\f{\bi{2k}k}k\eq-\l(\f{-1}p\r)\f83pE_{p-3}\pmod {p^2},$$
where $(\f{\cdot}p)$ denotes the Legendre symbol and $E_{p-3}$ stands for the $(p-3)$-th Euler number.

Recall that the Bernoulli polynomials are given by
$$B_n(x)=\sum_{k=0}^n\bi nkB_kx^{n-k}\ \ (n=0,1,2,\ldots).$$
Motivated by the above work, we mainly obtain the following result in this paper.

 \begin{theorem}\label{Th1.1} Let $p>3$ be a prime. Then
 \begin{equation}\label{1.1}\sum_{k=1}^{p-1}\f{\bi{2k}k}kH_k\eq\f13\l(\f p3\r)B_{p-2}\l(\f13\r)\pmod{p}
 \end{equation}
 and
 \begin{equation}\label{1.2}\sum_{k=1}^{p-1}\f{\bi{2k}k}kH_{2k}\eq\f7{12}\l(\f p3\r)B_{p-2}\l(\f13\r)\pmod{p}.
 \end{equation}
\end{theorem}

\begin{remark}\label{Rem1.1} Another motivation of Theorem 1.1 is our desire to prove Theorem 1.2 in this section.
Why (\ref{1.1}) and (\ref{1.2}) should involve the value of $B_{p-2}(x)$ at $1/3$ ? There are no intuitive reasons, but a slight indication
comes from the following congruence in [ST11, (1.5)] for any prime $p>3$:
$$\sum_{k=0}^{p-1}\bi {2k}k\eq\l(\f p3\r)+2p\sum_{k=1}^{p-1}\f{(-1)^k}k\l(\f{p-k}3\r)\pmod{p^3}.$$
There are no closed forms for the sums in (\ref{1.1}) and (\ref{1.2}).
Our approach to Theorem 1.1 is somewhat unique in the sense that it depends heavily on some special combinatorial identities.
\end{remark}

Clearly Theorem \ref{Th1.1} has the following consequence.

\begin{corollary}\label{Cor1.1} For any prime $p>3$ we have
\begin{equation}\label{1.3}\sum_{k=1}^{p-1}\f{\bi{2k}k}k(4H_{2k}-7H_k)\eq0\pmod p.
\end{equation}
\end{corollary}

Recall that $H_{p-1}\eq0\pmod{p^2}$ for any prime $p>3$. So our following conjecture is much stronger than Corollary \ref{Cor1.1}.
\begin{conjecture}\label{Conj1.1} For any prime $p>3$ we have
$$\sum_{k=1}^{p-1}\f{\bi{2k}k}k(4H_{2k}-7H_k)\eq-14\f{H_{p-1}}p+\f{278}{15}p^3B_{p-5}\pmod{p^4}.$$
\end{conjecture}

Note that Theorem \ref{Th1.1} and Corollary \ref{Cor1.1} are related to the second author's following conjectural formulas for $\zeta(3)=\sum_{n=1}^\infty1/n^3$
and $K=\sum_{k=1}^\infty(\f k3)/k^2$ (cf. \cite{S14a}):
$$\sum_{k=1}^\infty\f{H_{2k}+2H_k}{k^2\bi{2k}k}=\f 53\zeta(3)\ \t{and}\ \sum_{k=1}^\infty\f{H_{2k}+17H_k}{k^2\bi{2k}k}=\f 52\sqrt3\pi K.$$

The Franel numbers $f_n=\sum_{k=0}^n\bi nk^3\ (n=0,1,2,\ldots)$ play important roles in combinatorics and number theory.
In 1975 P. Barrucand \cite{B} obtained the identity $\sum_{k=0}^n\bi nk f_k=g_n$, where
\begin{equation}\label{1.4}g_n:=\sum_{k=0}^n\bi nk^2\bi{2k}k.
\end{equation}
The sequences $(f_n)_{n\gs0}$ and $(g_n)_{n\gs0}$ are two of the five sporadic sequences (cf. D. Zagier \cite[Section 4]{Z})
which are integral solutions of certain Ap\'ery-like recurrence equations and closely related to the theory of modular forms.
In 2013, Sun \cite{S13} revealed some unexpected connections between those numbers $f_n,g_n$ and representations of primes $p\eq1\pmod3$ in the form $x^2+3y^2$ with $x,y\in\Z$.
For any prime $p>3$, Sun \cite{S14b} and \cite[(1.15)]{S12b} showed that
$$\sum_{n=1}^{p-1}g_n\eq\sum_{n=1}^{p-1}h_n\eq0\pmod {p^2},$$
where
\begin{equation}\label{1.5}h_n:=\sum_{k=0}^n\bi nk^2C_k
\end{equation}
and $C_k$ refers to the Catalan number $\bi{2k}k/(k+1)=\bi{2k}k-\bi{2k}{k+1}$. The numbers $h_0,h_1,h_2,\ldots$ appeared naturally in the second author's study of Ap\'ery polynomials (cf. \cite{S12b}).

Applying Theorem \ref{Th1.1}, we deduce the following result.

\begin{theorem}\label{Th1.2} Let $p>3$ be a prime. Then
\begin{equation}\label{1.6}
\f1{p^2}\sum_{k=1}^{p-1}g_k\eq\sum_{k=1}^{p-1}g_kH_k^{(2)}\eq\frac{5}{8}\l(\frac{p}{3}\r)B_{p-2}\l(\frac{1}{3}\r)\pmod{p},
 \end{equation}
 and
\begin{equation}\label{1.7}
 \f1{p^2}\sum_{k=1}^{p-1}h_k\eq\sum_{k=1}^{p-1}h_kH_k^{(2)}\eq\frac{3}{4}\l(\frac{p}{3}\r)B_{p-2}\l(\frac{1}{3}\r)\pmod{p}.
 \end{equation}
 \end{theorem}

We are going to prove Theorems 1.1 and 1.2 in Sections 2 and 3 respectively. Our proofs make use of some sophisticated combinatorial identities.

 \section{Proof of Theorem 1.1}
 \setcounter{lemma}{0}
\setcounter{theorem}{0}
\setcounter{corollary}{0}
\setcounter{remark}{0}
\setcounter{equation}{0}
\setcounter{conjecture}{0}

 \begin{lemma}\label{Lem2.1} For any $n\in\N$, we have
 \begin{align}\label{2.1}\sum_{k=0}^n\bi xk\bi y{n-k}=&\bi{x+y}n,
\\\label{2.2} \sum_{k=0}^n \binom{n}{k}^2H_k=&\bi{2n}n(2H_n-H_{2n}),
\\\label{2.3} \sum_{k=0}^n(-1)^k\bi nk\bi{2k}k=&(-1)^n\sum_{k=0}^{\lfloor n/2 \rfloor}\bi n{2k}\bi{2k}k.
 \end{align}
\end{lemma}
 The three identities (\ref{2.1})-(\ref{2.3}) are known, see, e.g., \cite[(3.1), (3.125) and (3.86)]{G}.

 \medskip
\noindent{\it Proof of Theorem 1.1}. By (\ref{2.2}),
$$\sum_{j=1}^k\binom{k}{j}^2H_j=\bi{2k}k(2H_k-H_{2k})\quad\t{for each}\ k=1,\ldots,p-1.$$
Therefore
\begin{equation}\label{2.4}
\sum_{k=1}^{p-1}\f{\bi{2k}k}kH_{2k}=2\sum_{k=1}^{p-1}\f{\bi{2k}k}kH_k-\sum_{k=1}^{p-1}\f1k\sum_{j=1}^k\binom{k}{j}^2H_j.
\end{equation}
Observe that
$$\sum_{k=1}^{p-1}\f1k\sum_{j=1}^k\binom{k}{j}^2H_j=\sum_{j=1}^{p-1}\f{H_j}j\sum_{k=j}^{p-1}\bi kj\bi{k-1}{j-1}$$
and
\begin{align*}
&\sum_{k=j}^{p-1}\bi kj\bi{k-1}{j-1}
\\=&\sum_{i=0}^{p-1-j}\bi{i+j}i\bi{i+j-1}i=\sum_{i=0}^{p-1-j}\bi{-j-1}i\bi{-j}i
\\\eq&\sum_{i=0}^{p-1-j}\bi{p-j}i\bi{p-1-j}i=\sum_{i=0}^{p-1-j}\bi{p-j}i\bi{p-1-j}{p-1-j-i}
\\=&\bi{2p-2j-1}{p-1-j}\pmod{p}
\end{align*}
with the help of the Chu-Vandermonde identity (\ref{2.1}). Thus
\begin{align*}
\sum_{k=1}^{p-1}\f1k\sum_{j=1}^k\binom{k}{j}^2H_j
\eq&\sum_{j=1}^{p-1}\f{H_j}j\bi{2p-2j-1}{p-1-j}
\\\eq&\sum_{j=1}^{p-1}\f{H_j}j\bi{-2j-1}{p-1-j}
=\sum_{j=1}^{p-1}\f{H_j}j\bi{p+j-1}{2j}(-1)^j
\\=&\sum_{j=1}^{p-1}\f{H_j}j\cdot\f{p(-1)^j}{(2j)!(p+j)}\prod_{i=1}^j(p^2-i^2)
\\\eq&p\sum_{j=1}^{p-1}\f{H_j}{j^2\bi{2j}j}\eq p\sum_{j=(p+1)/2}^{p-1}\f{H_j}{j^2\bi{2j}j}\pmod p.
\end{align*}
By \cite[Lemma 2.1]{S11},
$$j\bi{2j}j\bi{2(p-j)}{p-j}\eq2p\pmod{p^2}\quad\t{for all}\ j=\f{p+1}2,\ldots,p-1.$$
(Tauraso \cite{T} contains a similar technique.)
Therefore
\begin{align*}&\sum_{k=1}^{p-1}\f1k\sum_{j=1}^k\binom{k}{j}^2H_j
\\\eq&\f 12\sum_{j=(p+1)/2}^{p-1}\f{H_j}j\bi{2(p-j)}{p-j}
=\f 12\sum_{k=1}^{(p-1)/2}\f{\bi{2k}kH_{p-k}}{p-k}\pmod p.
\end{align*}
Since
$$H_{p-k}=H_{p-1}-\sum_{0<j<k}\f1{p-j}\eq H_{k-1}=H_k-\f1k \pmod p$$
for all $k=1,\ldots,p-1$, from the above we obtain
\begin{align*}\sum_{k=1}^{p-1}\f1k\sum_{j=1}^k\binom{k}{j}^2H_j
\eq&-\f 12\sum_{k=1}^{p-1}
\f{\bi{2k}k}kH_k+\f12\sum_{k=1}^{p-1}\f{\bi{2k}k}{k^2}\pmod p.
\end{align*}
Combining this with (\ref{2.4}) we get
\begin{equation}\label{2.5}
\sum_{k=1}^{p-1}\f{\bi{2k}k}kH_{2k}\equiv\f 52\sum_{k=1}^{p-1}\f{\bi{2k}k}kH_k-\f 12\sum_{k=1}^{p-1}\f{\bi{2k}k}{k^2}\pmod p.
\end{equation}

For each $k=1,\ldots,p-1$, clearly
$$\bi pk=\f p k\prod_{0<j<k}\f{p-j}j\eq(-1)^{k-1}\f p k(1-pH_{k-1})\pmod{p^3}.$$
Thus
$$\sum_{k=1}^{p-1}(-1)^k\bi pk\bi{2k}k\eq-p\sum_{k=1}^{p-1}\f{\bi{2k}k}k(1-pH_{k-1})\pmod{p^3}.$$
On the other hand, by (\ref{2.3}) we have
\begin{align*}\sum_{k=0}^p(-1)^k\bi pk\bi{2k}k=&(-1)^p\sum_{k=0}^{(p-1)/2}\bi p{2k}\bi{2k}k
\\\eq&-1+p\sum_{k=1}^{(p-1)/2}\f{1-pH_{2k-1}}{2k}\bi{2k}k\pmod {p^3}.
\end{align*}
Therefore
\begin{equation}\label{2.6}\begin{aligned}&-p\sum_{k=1}^{p-1}\f{\bi{2k}k}k\l(1-pH_k+\f pk\r)-\bi{2p}p+1
\\\eq&-1+p\sum_{k=1}^{(p-1)/2}\f{1-p(H_{2k}-1/(2k))}{2k}\bi{2k}k\pmod {p^3}
\end{aligned}\end{equation}
Since $\bi{2p}p\eq2\pmod{p^3}$ by Wolstenholme's theorem, and
\begin{equation}\label{2.7}\sum_{k=1}^{p-1}\f{\bi{2k}k}k\eq0\pmod{p^2}
\end{equation}
by \cite{ST10}, from (\ref{2.6}) we get
\begin{equation}\label{2.8}
\sum_{k=1}^{p-1}\frac{\binom{2k}{k}}{k}H_k\equiv\frac{1}{2p}\sum_{k=1}^{(p-1)/2}\frac{\binom{2k}{k}}{k}-\frac{1}{2}\sum_{k=1}^{(p-1)/2}\frac{\binom{2k}{k}}{k}H_{2k}+
\frac{5}{4}\sum_{k=1}^{(p-1)/2}\frac{\binom{2k}{k}}{k^2}\pmod p.
\end{equation}
(Note that $\sum_{k=1}^{p-1}\bi{2k}k/k^2\eq\sum_{k=1}^{(p-1)/2}\bi{2k}k/k^2\pmod p$.).

Clearly,
$$pH_{2p-2k}=p\sum_{j=1\atop j\not=p}^{2p-2k}\f1j+1\eq1\pmod{p}$$
for all $k=1,\ldots,(p-1)/2$. So we have
\begin{align*}
\sum_{k=1}^{(p-1)/2}\f1{k^2\bi{2k}k}\eq&\sum_{k=1}^{(p-1)/2}\frac{pH_{2p-2k}}{k^2\binom{2k}{k}}
=\sum_{j=(p+1)/2}^{p-1}\frac{pH_{2j}}{(p-j)^2\binom{2(p-j)}{p-j}}
\\\eq&\f1{2}\sum_{j=(p+1)/2}^{p-1}\f{\bi{2j}j}jH_{2j}\pmod p
\end{align*}
with the help of \cite[Lemma 2.1]{S11}.
By \cite[(1.2) and (1.3)]{S11},
$$\frac{1}{p}\sum_{k=1}^{(p-1)/2}\frac{\binom{2k}{k}}{k}+\sum_{k=1}^{(p-1)/2}\frac{2}{k^2\binom{2k}{k}}\equiv0\pmod p.$$
Therefore
$$\sum_{k=1}^{p-1}\frac{\binom{2k}{k}}{k}H_{2k}\equiv\sum_{k=1}^{(p-1)/2}\frac{\binom{2k}{k}}{k}H_{2k}-\frac{1}{p}\sum_{k=1}^{(p-1)/2}\frac{\binom{2k}{k}}{k}\pmod p.$$
Combining this with (\ref{2.8}) we get
\begin{equation}\label{2.9}
\sum_{k=1}^{p-1}\frac{\binom{2k}{k}}{k}H_k\equiv\frac{5}{4}\sum_{k=1}^{p-1}\frac{\binom{2k}{k}}{k^2}-\frac{1}{2}\sum_{k=1}^{p-1}\frac{\binom{2k}{k}}{k}H_{2k}\pmod p.
\end{equation}

(\ref{2.5}) and (\ref{2.9}) together imply that
$$\sum_{k=1}^{p-1}\f{\bi{2k}k}kH_k\eq\f23\sum_{k=1}^{p-1}\f{\bi{2k}k}{k^2}\pmod p
\ \t{and}\ \sum_{k=1}^{p-1}\f{\bi{2k}k}kH_{2k}\eq\f76\sum_{k=1}^{p-1}\f{\bi{2k}k}{k^2}\pmod p.$$
It is known that
\begin{equation}\label{2.10}\sum_{k=1}^{p-1}\f{\bi{2k}k}{k^2}\eq\f12\l(\f p3\r)B_{p-2}\l(\f13\r)\pmod p
\end{equation}
(cf. \cite{MT}). So we get the desired (\ref{1.1}) and (\ref{1.2}).

\begin{remark}\label{Rem2.1} In \cite{S11} the second author proved that
$$-\f1p\sum_{k=1}^{(p-1)/2}\f{\bi{2k}k}k\eq\sum_{k=1}^{(p-1)/2}\f2{k^2\bi{2k}k}\eq\l(\f{-1}p\r)\f83E_{p-3}\pmod p.$$

\end{remark}
 \section{Proof of Theorem 1.2}
 \setcounter{lemma}{0}
\setcounter{theorem}{0}
\setcounter{corollary}{0}
\setcounter{remark}{0}
\setcounter{equation}{0}
\setcounter{conjecture}{0}

 \begin{lemma}\label{Lem3.1} For any nonnegative integers $m$ and $n$, we have
 \begin{equation}\label{3.1}
 \sum_{k=0}^n\binom{x+k}{m}=\binom{n+x+1}{m+1}-\binom{x}{m+1}.
 \end{equation}
 and
 \begin{equation}\label{3.2}\sum_{k=0}^n\binom{n}{k}^2\binom{x+k}{2n}=\binom{x}{n}^2.
 \end{equation}

 \end{lemma}
\begin{remark}\label{Rem3.1} Both (\ref{3.1}) and (\ref{3.2}) can be found in \cite[(1.48) and (6.30)]{G}.
\end{remark}

 \begin{lemma}\label{Lem3.2} For any nonnegative integer $n$, we have
\begin{equation}\label{3.3}
 \sum_{k=0}^n\binom{n}{k}^2\binom{x+k}{2n+1}=\f1{(4n+2)\bi{2n}n}\sum_{k=0}^n(2x-3k)\bi xk^2\bi{2k}k.
 \end{equation}
 \end{lemma}
 \begin{remark}\label{Rem3.2} One might wonder how we find (\ref{3.3}). In fact, we use the software package {\tt Sigma} to find
 a recurrence for $u_n=\sum_{k=0}^n\binom{n}{k}^2\binom{x+k}{2n+1}$ and then obtain (\ref{3.3}) by solving the recurrence for $u_n$ via {\tt Sigma}.
 \end{remark}

 \noindent{\it Proof of Lemma 3.2}. Let $F(x)$ and $G(x)$ denote the left-hand side and the right-hand side of (\ref{3.3}). With the help of (\ref{3.2}), we see that
 $$F(x+1)-F(x)=\bi xn^2.$$
 Applying the Zeilberger algorithm (see, e.g., \cite[pp.\,101-119]{PWZ}) via {\tt Mathematica}, we find that
 $$G(x+1)-G(x)=\bi xn^2\quad\t{for all}\ x=0,1,2,\ldots.$$
 So, by induction $F(x)=G(x)$ for all $x\in\N$. As $F(x)$ and $G(x)$ are polynomials in $x$ of degree $2n+1$, we have the desired (\ref{3.3}).
 \qed

 \medskip\noindent{\it Proof of Theorem 1.2}. (i) With the help of Lemma \ref{Lem3.1}, we have
 \begin{align*}
 \sum_{n=0}^{p-1}g_n=&\sum_{n=0}^{p-1}\sum_{k=0}^n\binom{n}{k}^2\binom{2k}{k}=\sum_{k=0}^{p-1}\binom{2k}{k}\sum_{n=k}^{p-1}\binom{n}{k}^2
 \\=&\sum_{k=0}^{p-1}\binom{2k}{k}\sum_{n=k}^{p-1}\sum_{j=0}^k\binom{k}{j}^2\binom{n+j}{2k}
 \\=&\sum_{k=0}^{p-1}\binom{2k}{k}\sum_{j=0}^k\binom{k}{j}^2\sum_{n=k}^{p-1}\binom{n+j}{2k}
 \\=&\sum_{k=0}^{p-1}\binom{2k}{k}\sum_{j=0}^k\binom{k}{j}^2\binom{p+j}{2k+1}.
 \end{align*}
Thus, by applying Lemma \ref{Lem3.2} we get

     \begin{align*}\sum_{k=0}^{p-1}g_k=&\sum_{k=0}^{p-1}\f1{4k+2}\sum_{j=0}^k(2p-3j)\bi pj^2\bi{2j}j
     \\=&\f12\sum_{j=0}^{p-1}(2p-3j)\bi pj^2\bi{2j}j\(\sum_{k=0}^{p-1}\f1{2k+1}-\sum_{0\ls i<j}\f1{2i+1}\)
     \\=&\f12\sum_{j=0}^{p-1}(2p-3j)\bi pj^2\bi{2j}j\l(H_{2p-1}-\f{H_{p-1}}2-H_{2j}+\f{H_j}2\r).
     \end{align*}
     Note that $pH_{p-1}\eq0\pmod{p^3}$ and
     \begin{align*}pH_{2p-1}=&1+p\sum_{j=1}^{p-1}\l(\f1{p-j}+\f1{p+j}\r)=1+p\sum_{j=1}^{p-1}\f{2p}{p^2-j^2}
     \\\eq&1-2p^2H_{p-1}^{(2)}\eq1\pmod{p^3}.
     \end{align*}
     Therefore
\begin{align*}\sum_{k=0}^{p-1}g_k\eq&\sum_{j=0}^{p-1}\f{2p-3j}{2p}\bi pj^2\bi{2j}j
\\&+\sum_{j=0}^{p-1}\f{2p-3j}2\bi pj^2\bi{2j}j\l(\f{H_j}2-H_{2j}\r)
\pmod{p^3}
\end{align*}
and hence
\begin{align*}\sum_{k=1}^{p-1}g_k\eq&\sum_{j=1}^{p-1}\f{2p-3j}{2p}\cdot\f{p^2}{j^2}\bi {p-1}{j-1}^2\bi{2j}j
\\&+\sum_{j=1}^{p-1}\f{2p-3j}2\cdot\f{p^2}{j^2}\bi{p-1}{j-1}^2\bi{2j}j\l(\f{H_j}2-H_{2j}\r)
\\\eq&p^2\sum_{j=1}^{p-1}\f{\bi{2j}j}{j^2}-\f{3p}2\sum_{j=1}^{p-1}\f{\bi{2j}j}j\bi{p-1}{j-1}^2
\\&-\f32p^2\sum_{j=1}^{p-1}\f{\bi{2j}j}j\l(\f{H_j}2-H_{2j}\r)
\pmod{p^3}.
\end{align*}
(Note that $\bi{2j}jH_{2j}$ is $p$-adic integral for all $j=1,\ldots,p-1$.) Clearly,
\begin{equation}\label{3.4}\bi{p-1}{j-1}^2\eq(1-pH_{j-1})^2\eq1-2pH_{j-1}\pmod{p^2}.
\end{equation}
Thus
\begin{align*}\f1{p^2}\sum_{k=1}^{p-1}g_k\eq&\sum_{k=1}^{p-1}\f{\bi{2k}k}{k^2}-\f32\(\f1p\sum_{k=1}^{p-1}\f{\bi{2k}k}k-2\sum_{k=1}^{p-1}\f{\bi{2k}k}k\l(H_k-\f1k\r)\)
\\&-\f32\(\f12\sum_{k=1}^{p-1}\f{\bi{2k}k}kH_k-\sum_{k=1}^{p-1}\f{\bi{2k}k}kH_{2k}\)
\\\eq&-2\sum_{k=1}^{p-1}\f{\bi{2k}k}{k^2}+\f 94\sum_{k=1}^{p-1}\f{\bi{2k}k}kH_k+\f32\sum_{k=1}^{p-1}\f{\bi{2k}k}kH_{2k}\pmod p
\end{align*}
with the help of (\ref{2.7}). Now, applying Theorem \ref{Th1.1} and (\ref{2.10}) we immediately get that
\begin{equation}\label{3.5}\f1{p^2}\sum_{k=1}^{p-1}g_k\eq\f58\l(\f p3\r)B_{p-2}\l(\f13\r)\pmod p.
\end{equation}

(ii) Observe that
\begin{align*}\sum_{n=0}^{p-1}(2g_n-h_n)=&\sum_{n=0}^{p-1}\sum_{k=0}^n\bi nk^2\l(2-\f1{k+1}\r)\bi{2k}k
\\=&\sum_{k=0}^{p-1}\f{2k+1}{k+1}\bi{2k}k\sum_{n=k}^{p-1}\bi nk^2.
\end{align*}
Similar to the proof in part (i), we have
\begin{align*}\sum_{n=0}^{p-1}(2g_n-h_n)=&\sum_{k=0}^{p-1}\f{(2k+1)/(k+1)}{4k+2}\sum_{j=0}^k(2p-3j)\bi pj^2\bi{2j}j
\\=&\f12\sum_{j=0}^{p-1}(2p-3j)\bi pj^2\bi{2j}j\l(H_{p-1}+\f1p-H_j\r)
\end{align*}
and thus
\begin{align*}&\sum_{k=1}^{p-1}(2g_k-h_k)
\\\eq&\f12\sum_{j=1}^{p-1}(2p-3j)\f{p^2}{j^2}\bi{p-1}{j-1}^2\bi{2j}j\l(\f1p-H_j\r)
\\\eq&p^2\sum_{j=1}^{p-1}\f{\bi{2j}j}{j^2}-\f{3p}2\sum_{j=1}^{p-1}\f{\bi{2j}j}j\l(1-2pH_{j-1}\r)(1-pH_j)\ (\t{by (\ref{3.4})})
\\\eq&p^2\sum_{j=1}^{p-1}\f{\bi{2j}j}{j^2}-\f32p\sum_{j=1}^{p-1}\f{\bi{2j}j}j\l(1+p\l(\f2j-3H_j\r)\r)
\\=&-2p^2\sum_{j=1}^{p-1}\f{\bi{2j}j}{j^2}-\f 32p\sum_{j=1}^{p-1}\f{\bi{2j}j}j+\f 92p^2\sum_{j=1}^{p-1}\f{\bi{2j}j}jH_j\pmod{p^3}
\end{align*}
Combining this with (\ref{1.1}), (\ref{2.7}) and (\ref{2.10}), we obtain that
$$\sum_{k=1}^{p-1}(2g_k-h_k)\eq\f{p^2}2\l(\f p3\r)B_{p-2}\l(\f13\r)\pmod{p^3}.$$
Thus,
\begin{equation}\label{3.6}\sum_{k=1}^{p-1}h_k\eq\f 34p^2\l(\f{p}3\r)B_{p-2}\l(\f13\r)\pmod{p^3}
\end{equation}
with the help of (\ref{3.5}).

(iii) By \cite[Theorem1.1]{S14b},
\begin{equation}\label{3.7}\sum_{k=1}^{p-1}g_k\eq p^2\sum_{k=1}^{p-1}g_kH_k^{(2)}+\frac{7}{6}p^3B_{p-3}\pmod{p^4}.
\end{equation}
Therefore
$$\sum_{k=1}^{p-1}g_kH_k^{(2)}\eq\f1{p^2}\sum_{k=1}^{p-1}g_k\pmod p$$
and hence (\ref{1.6}) holds in view of (\ref{3.5}).

 From \cite[Theorem1.1]{S14b} we know that
\begin{equation}\label{3.8}\sum_{k=0}^{p-1}g_k(x)\l(1-p^2H_k^{(2)}\r)\eq\sum_{k=0}^{p-1}\f p{2k+1}\(1-2p^2H_k^{(2)}\)x^k\pmod{p^4},
 \end{equation}
 where $g_k(x)=\sum_{j=0}^k\bi{k}{j}^2\bi{2j}jx^j$.
 Therefore, the left-hand side of (\ref{3.8}) minus the right-hand side of (\ref{3.8}) can be written as $p^4P(x)$, where $P(x)$
 is a polynomial of degree at most $p-1$ with $p$-adic integer coefficients.
Since
     $$h_n=\sum_{k=0}^n\binom{n}{k}^2\bi{2k}kC_k=\int_0^1g_n(x)dx\quad\t{for}\ n=0,1,2,\ldots,$$
we deduce that
     \begin{align*}
     &\sum_{k=0}^{p-1}h_k(1-p^2H_k^{(2)})
     \\=&\int_0^1\sum_{k=0}^{p-1}g_k(x)(1-p^2H_k^{(2)})dx
     \\=&\sum_{k=0}^{p-1}\f p{2k+1}(1-2p^2H_k^{(2)})\int_0^1x^kdx+p^4\int_0^1P(x)dx
     \\\eq&\sum_{k=0}^{p-1}\f{2p}{2k+1}(1-2p^2H_k^{(2)})-\sum_{k=0}^{p-1}\f p{k+1}(1-2p^2H_k^{(2)})\pmod{p^3}.
     \end{align*}
 Combining (\ref{3.7}) and (\ref{3.8}) we see that
 \begin{align*}1+\f 76p^3B_{p-3}\eq&\sum_{k=0}^{p-1}g_k\l(1-p^2H_k^{(2)}\r)
 \\\eq&\sum_{k=0}^{p-1}\f p{2k+1}\l(1-2p^2H_k^{(2)}\r)\pmod{p^4}.
 \end{align*}
 Therefore
\begin{align*}
\sum_{k=0}^{p-1}h_k(1-p^2H_k^{(2)})\eq&2+\f73p^3B_{p-3}-\sum_{k=0}^{p-1}\f p{k+1}+2p^3\sum_{k=0}^{p-1}\f{H_k^{(2)}}{k+1}
     \\\eq&2-1-pH_{p-1}+2p^2H_{p-1}^{(2)}\eq1\pmod{p^3}
     \end{align*}
which implies that $$\sum_{k=1}^{p-1}h_k(1-p^2H_k^{(2)})\eq 0\pmod{p^3}.$$
Combining this (\ref{3.6}) we obtain the desired (\ref{1.7}).

So far we have completed the proof of Theorem \ref{Th1.2}. \qed

\Ack. The authors would like to thank the referee for helpful comments.

     \end{document}